\documentclass[11pt,twoside]{amsart}
\usepackage{amsmath, amsthm, amsfonts, amssymb, amscd}
\usepackage{graphicx}
\usepackage[ruled,vlined,norelsize]{algorithm2e}
\usepackage{caption}
\usepackage{subcaption}
\usepackage{mathtools}
\usepackage{tikz}
\usepackage{hyperref}
\usepackage{comment}
\usepackage{xcolor}

\theoremstyle{plain}
\newtheorem{theorem}{Theorem}[section]

\theoremstyle{definition}

\theoremstyle{remark}

\newcommand{\field}[1]{\mathbb{#1}}
\newcommand{\R}{\field{R}}

\def \bv {{\bf v}}

\def \bx {{\bf x}}

\def \bF {{\bf F}}

\def \cD {\mathcal{D}}
\def \cE {\mathcal{E}}

\def \cK {\mathcal{K}}

\def \cP {\mathcal{P}}

\def \dt  {\, \mbox{d}t}

\def \ddt  {\frac{\mbox{d\,\,}}{\mbox{d}t}}

\def \d {\delta}

\def \e {\varepsilon}

\def \vmax {\bv_{\mathrm{max}}}
\def \vmin {\bv_{\mathrm{min}}}
\def \Fmax {\bF_{\mathrm{max}}}

\title{Collective Behavior in Systems within Confined Environments}
\author{Veronica Kalicki}
\address{Veronica Kalicki, Departments of Mathematical Sciences and Compute Science, University of Illinois at Chicago}
\date{May 1, 2020}

\begin{document} 

\begin{abstract}
In this paper, we study the behavior of systems of individuals in confined environments that are driven by laws of self-organization. We propose that, under certain conditions, the long-term behavior of such systems will be global alignment. We study the result by Felipe Cucker and Steve Smale as well as their models describing the evolution of a flock in continuous and discrete time. Specifically, we will describe the models of Cucker-Smale on bounded domains $\mathbb{R}^+$ and on the interval $[a,b]$. 
\end{abstract}

\maketitle
\tableofcontents
\section{Introduction}
Consider a system of individual ``agents" enclosed in a confined environment, driven by laws of self-organization. One may ask whether mathematical emergence may allow for a purely local communication protocol between agents to result in global collective behavior. The Cucker-Smale system introduced in \cite{CS2007a} and \cite{CS2007b} show that such systems, in the setting of Euclidean space, exhibit long-term global alignment. In this project, we adapt their methods to the settings of the half line and closed interval, which notably have nonempty boundary.
\\
\\

\section{One wall dynamics: flocks on $\R^+$}  
Let us consider Cucker-Smale system on the half-line $\mathbb{R}^+$: 

\begin{equation}\label{Cucker-Smale}
\begin{cases} 
      \dot{\bx}_i&=\bv_i,\\
      \dot{\bv}_i&= \frac{1}{N} \sum_{j=1}^N \phi(\bx_i-\bx_j)(\bv_j-\bv_i) +\bF_i
   \end{cases}
\end{equation} 

where $\bF_i$ is a potential force that repels the agents from the wall $x=0$. We consider it to be given by a potential $U$: 
\[ -U'(\bx_i) = \bF_i,\]
where $U$ is a decreasing function with 
\[ \sup U\subset [0,\ell), \indent \lim_{x\to 0} U(\bx)=\infty. \]
So, $\bF_i$ is a positive force which intensifies as $x\to 0$. Thus it ``pushes'' agents away from the wall. The parameter $\ell$ defines a reaction length scale and is a property of agents themselves. We assume the kernel $\phi$ is smooth and only depends on the absolute value of the input. We also assume that the system has a ``fat tail'', i.e. 
\begin{equation}\label{Fat Tail}
\int_0^\infty \phi(x) dx=\infty
\end{equation}  the so-called ``fat tail'' condition. \\

\indent The full energy of the system is given by, 
\begin{equation}\label{Energy Key}
\cE= \cK+ \cP
\end{equation}

\[\cK= \frac{1}{2N} \sum_{i=1}^N |\bv_i|^2, \indent \cP= \frac{1}{N} \sum_{i=1}^N U(\bx_i).\]

Another important quantity we will track is the total momentum of the system:
\begin{equation}\label{momentum}
p=\frac{1}{N} \sum_{i=1}^N \bv_i.
\end{equation}

\begin{theorem}
Any solution to the system \eqref{Cucker-Smale} with initial condition $\bx_i(0) >0$ will fulfill the following dynamic behavior:
\begin{itemize}
	\item[(1)] It will never collide with the wall, $\bx_i(t) >0$ for all $i=1,...,N$ and all $t >0$.
	\item[(2)] It will align
	\[
	A(t) = \max_{i,j} |\bv_i - \bv_j| \to 0,
	\]
	and flock strongly
	\[
	\bx_i(t) - \bx_j(t) \to \overline{ \bx}_{ij}(t)
	\] 
	with all agents settling outside the range of influence of the wall
	\[
	\liminf_{t \to \infty} \bx_i(t) \geq \ell.
	\]
	\item[(3)] Moreover, if the initial momentum points away from the wall, $p_0 >0$, then the flock will escape the influence of the wall in finite time, and the alignment will take place exponentially fast:
	\[
	A(t) \leq C e^{-\d t}.
	\]
\end{itemize}
\end{theorem}

The remainder of this section is dedicated to the proof of this theorem.

 We derive the equation for energy $\ddt \cE$.
Since,  \[ \ddt \cE= \ddt(\cK+\cP) = \ddt \cK + \ddt \cP \] 
and, \[ \ddt \cK= \ddt \big(\frac{1}{2N} \sum_{i=1}^N |\bv_i|^2\big),\indent  \ddt \cP=\ddt \big(\frac{1}{N} \sum_{i=1}^N U(\bx_i)\big). \] 

We obtain, 
\[ \begin{array}{lcl}
\ddt \cE & = & \ddt \cK + \ddt \cP \\
&=& \displaystyle \frac{1}{2N} \sum_{i=1}^N 2 {\bv_i}{\dot{\bv_i}} + \frac{1}{N} \sum_{i=1}^N \dot{U}(\bx_i)\bv_i \\
&\\
&=& \displaystyle \frac{1}{N} \sum_{i=1}^N \left(\dot{\bv_i}+ \dot{U}(\bx_i) \right) \bv_i \\
&=& \displaystyle \frac{1}{N} \sum_{i=1}^N (\dot{\bv_i}- \bF_i) v_i \\
&=& \displaystyle \frac{1}{N} \sum_{i=1}^N \left( \frac{1}{N} \sum_{j=1}^N \phi(\bx_i-\bx_j)(\bv_j-\bv_i) \right) \bv_i\\
&=& \displaystyle \frac{1}{N^2} \sum_{i,j=1}^N \phi(\bx_i-\bx_j)(\bv_j-\bv_i) \bv_i 
\end{array} \]
Now, we write \\
\[ \begin{array}{lcl}
(\bv_j- \bv_i)\bv_i & = & (\bv_j-\bv_i)(\bv_i-\bv_j) +(\bv_j-\bv_i) \bv_j \\
&= & -|\bv_j -\bv_i|^2- (\bv_i -\bv_j) \bv_j. 
\end{array}\]
Thus, 
\[ \begin{array}{lcl}
\ddt \cE &=& \displaystyle  \frac{1}{N^2} \sum_{i,j=1}^N \phi(\bx_i-\bx_j)(\bv_j-\bv_i)\bv_i \\
&=& \displaystyle \frac{-1}{N^2} \sum_{i,j=1}^N \phi(\bx_i-\bx_j)|\bv_i-\bv_j|^2 - \frac{1}{N^2} \sum_{i,j=1}^N \phi(\bx_i - \bx_j)(\bv_i-\bv_j)\bv_j 
\end{array} \]
The second sum is again $\ddt \cE,$ thus:
\[ \ddt \cE = \frac{-1}{2N^2}\sum_{i,j=1}^N \phi(\bx_i-\bx_j)|\bv_i-\bv_j|^2 \coloneqq -I_2.\]
\\
\\
We now show that the agents of the system never hit the wall. Recall that 
$\cP = \displaystyle \frac{1}{N} \sum_{i=1}^N U(\bx_i) $. In particular, $\cP \ge \frac{1}{N} U(\bx_i)$ for any i. Since $\lim_{x\to 0} U(\bx) = \infty$, we may choose some $\e$ such that $|\bx| < \e$ implies $U(\bx) \ge N(G+1)$, where $G$ is the initial energy. If the ith agent hits a wall, then $\bx_i(t) < \e$ for some $t$. This would imply $U(\bx_i(t)) \ge N(G+1)$, and thus that $\cP \ge G+1$. However, the total energy of the system is decreasing by the previous computation. Thus, the agents can not hit the wall, so the system has global solutions for any initial data. \\
\\
Again because the energy is decreasing we know that: 
$$\frac{1}{2N}\sum_{i=1}^N |\bv_i|^2 = K\leq G.$$ This implies in particular that the velocities are uniformly bounded by $\sqrt{2NG}$. \\
\\
We show that the diameter of the flock grows at most linearly in time:
\[ \begin{array}{lcl}
\cD(t) &= & \displaystyle \max_{i,j} | \bx_i-\bx_j| \le \max_{i,j} |\bx_i(t)|+ |\bx_j(t)| \\
\dot{\cD}(t) & \le & \displaystyle \max_{i,j} \frac{d}{dt} |\bx_i(t)|+ \frac{d}{dt}|\bx_j(t)| \le 2\sqrt{2NG} \\
\cD(t) & \le &  \displaystyle 2\sqrt{2NG}\thinspace t + D_0.
\end{array} \] \\
\\
Let us consider the evolution of total momentum now:
\[ \ddt p= \frac{1}{N} \sum_{i=1}^N \bF_i:=\bF \]
Since all the forces are positive, the momentum is increasing as long as there are agents in the $\ell$-vicinity of the wall. Moreover, integrating the momentum equation we obtain
\[ p(t)-p_0 = \int_0^t \bF(s) ds.\]
Since the velocities are uniformly bounded and hence so is momentum, we obtain global integrability of the force 
\begin{equation}\label{global integrability of force}
\int_0^\infty \bF(s) ds < \infty.
\end{equation}
With this information in mind we make maximum principle computation. We denote
\[
\vmax = \max_i \bv_i, \quad \vmin = \min_i \bv_i, \quad \Fmax =  \max_i \bF_i
\]
We now find bounds on the derivatives of the maximum and minimum velocities. For the maximum velocity, we have:
\[\ddt \vmax  \leq \phi(\cD) ( p - \vmax) + \Fmax\]
To prove this, first note that: \\
 \[\dot{\bv_i} = \frac{1}{N} \sum_{j=1}^N \phi(\bx_i-\bx_j)\cdot (\bv_j-\bv_i)+\bF_i. \]
Assume, $\dot{\bv_i}=\vmax$ \\

\[ \begin{array}{lcl}
\ddt{\vmax} &=& \displaystyle \frac{1}{N} \sum_{j=1}^N \phi(\bx_i-\bx_j)\cdot (\bv_j-\vmax)+ \bF_i \\
& \leq & \displaystyle \frac{1}{N}\sum_{j=1}^N \phi(\bx_i-\bx_j)\cdot (\bv_j- \vmax)+ \Fmax \\
& \leq & \displaystyle \frac{1}{N} \sum_{j=1}^N \phi(\cD(t))\cdot (\bv_j- \vmax) + \Fmax \\
& =& \displaystyle \phi(\cD)\cdot \frac{1}{N}( \sum_{j=1}^N \bv_j-\vmax) + \Fmax\\
 & =&\displaystyle \phi(\cD)\cdot(p(t)- \vmax)+\Fmax 
\end{array}\]\\

Let $\dot{\bv_i}=\vmin$, and note that $(\bv_j- \vmin)$ is positive

\[ \begin{array}{lcl}
\ddt \vmin & = & \displaystyle \frac{1}{N} \sum_{j=1}^N \phi(\bx_i-\bx_j)\cdot (\bv_j-\vmin)+\bF_i \\ 
&\geq & \displaystyle \frac{1}{N} \sum_{j=1}^N \phi(\cD)\cdot (\bv_j-\vmin) .
\end{array} \]

Since $\phi$ is decreasing, its minimum is achieved at $\cD=\max|\bx_i-\bx_j|$. \\
\\
Let us consider the following amplitude, 

\[\begin{array}{lcl}
A& =& \displaystyle \vmax - \vmin \\
\ddt A& =& \displaystyle \ddt \vmax - \ddt \vmin  \\
& \leq & \displaystyle  \phi(\cD)\cdot (p-\vmax)+ \Fmax -\phi(\cD) \cdot (p-\vmin) \\
&= & \displaystyle \phi(\cD)\cdot(\vmin-\vmax)+ \Fmax. 
\end{array} \]\\
\\
Here, $(\vmin-\vmax)=-A$. Thus we obtain, $\ddt A\leq -\phi(\cD)A+ \Fmax$ where the force $\Fmax$ is integrable. \\
\\
\\
Recall that $ \displaystyle \ddt A \le -\phi(\cD)A + \Fmax$ and $\ddt \cD \le A.$ Now consider the equation:
\begin{equation}\label{L}
 \displaystyle L(A(t),\cD(t)) = A(t) + \int_0^D(t) \phi(r)dr.
\end{equation} Differentiating each side and applying the above inequalities, we obtain:\\
\[ \begin{array}{lcl}
\ddt L& =& \displaystyle \ddt A + \phi(D)A \\
&\le& \displaystyle -\phi(D)A + \Fmax + \phi(D)A  \\
&=& \displaystyle \Fmax 
\end{array} \]
\\
Thus, $ \displaystyle L(A(t),\cD(t)) \le L(A(0),\cD(0)) + \int_0^t \Fmax(t)dt.$ Using the previous expression for $L(A(t),D(t))$, we obtain: \\
\\
\[ \displaystyle A(t) + \int_0^{D(t)} \phi(r)dr \le A(0) + \int_0^{D(0)} \phi(r)dr + \int_0^t \Fmax(t)dt\] 
\\
In particular, since $A(t)$ is bounded and the total force is integrable, we have that $\displaystyle\int_0^{\cD(t)} \phi(r)dr$ is bounded. By the fat tail condition, this implies that $\cD(t)$ must also be bounded, which is what it means for the system to exhibit flocking. In what follows, we let $\overline{\cD}$ be an upper bound for $\cD(t).$ \\
\\
Since the system flocks, we know that $\ddt A \le -cA + \Fmax$ for some fixed constant $c$. Letting $\mathbb{F} = \displaystyle \int_0^{\infty} \Fmax(t) dt < \infty$ and integrating both sides of the inequality, we obtain:\\
\[  A(t) \le A_0 e^{-ct} + \int_0^t e^{-c(t-s)} F(s)ds. \]
\\

It suffices to show the integral converges to 0. 
We split it up as:
\[\displaystyle \int_0^{t/2} e^{-c(t-s)} F(s)ds + \int_{t/2}^t e^{-c(t-s)} F(s)ds.\]
\\ 
For $0 \le s \le t/2$, we have $t-s \ge t/2$. Thus, the first integral is bounded by $\displaystyle e^{-ct/2}\int_0^{t/2} F(s)ds \le e^{-ct/2}\mathbb{F}$, which goes to $0$ as $t$ goes to infinity. Now, note that $e^{-c(t-s)} \le 1$ for all $s \le t$, so the second integral is bounded by $\displaystyle\int_{t/2}^t F(s)ds$, which also goes to $0$ as $t$ goes to infinity. Thus, so does $A(t)$.\\
\\
\\
Now, we wish to show the system exhibits strong flocking, i.e. the amplitude is globally integrable. By the previous discussion, we have
\[ \int_0^{\infty} A(t) dt \le \frac{1}{c}A_0 + \int_0^{\infty} \int_0^t e^{-c(t-s)} ds dt .\]
By Fubini's theorem, we may rewrite the latter integral as:
\[ \begin{array}{lcl}
\displaystyle \int_0^{\infty} \int_s^{\infty} F(s)e^{-c(t-s)} dt ds& =& \displaystyle \int_0^{\infty} F(s) \int_s^{\infty} e^{-c(t-s)} \\
&=& \displaystyle \frac{1}{c} \int_0^{\infty} F(s) ds \\
&\le& \displaystyle \frac{1}{c}\mathbb{F}
\end{array} \]
Thus, the system exhibits strong flocking. \\
\\
Now, we wish to show that flock stabilizes, i.e. that the distances between the agents converge:
\begin{equation}\label{agents converge}
\bx_i(t) - \bx_j(t) \to \overline{ \bx}_{ij}(t).
\end{equation}
To do this, write:
\[\bx_i(t) - \bx_j(t) = \bx_i(0) - \bx_j(0) + \int_0^t \bv_i(s) - \bv_j(s) ds.\]
This integral is bounded in absolute value by the amplitude, which we previously established is globally integrable. Thus, the flock stabilizes.\\
\\
If we suppose the initial momentum is positive $p_0 > 0$, then ``on average" the flock is heading to the right away from the wall. We would expect that in this case eventually the flock would escape the influence of the wall and from some finite time $t^*$ will evolve as in the open space (let's call it ``free flock").  Since the momentum is increasing, $p(t) \geq p_0$ for all $t>0$. Then from the above ,
\[
\ddt \vmin  \geq \phi(\bar{\cD}) ( p_0 - \vmin) .
\]
We now show, using Duhamel's principle, that there exists some finite $t^*$ such that $\vmin \ge \frac{1}{2} p_0$ for $t > t^*$.\\
\\
Take, \[\bv(t)\geq e^{-c_0t}\bv_0+ \int_0^te^{-c_0(t-s)} c_0\thinspace p_0 \thinspace ds = e^{-c_0t}\bv_0+ \frac{c_0p_0}{c_0}(1-e^{-c_0t}). \]
Note that as $t$ tends to infinity, \[e^{-c_0t}\bv_0\to 0 \indent \text{and} \indent \frac{c_0p_0}{c_0}(1-e^{-c_0t})\to p_0.\] \\
Since \[e^{-c_0t} \bv_0 + p_0(1-e^{-c_0t})\] converges to $p_0$, it is after some finite time $t_{\e}$ bounded below by $p_0 - \e$. Then we may let $t^* = t_{1/2}$.\\
\\
This will prove that all agents $\bx_i$ will have positive velocities eventually and hence will escape the interval $[0,\ell)$. From that time on the equation becomes the classical forceless Cucker-Smale:
\[ \begin{cases} 
      \dot{\bx}_i&=\bv_i,\\
      \dot{\bv}_i&= \frac{1}{N} \sum_{j=1}^N \phi(\bx_i-\bx_j)(\bv_j-\bv_i) 
   \end{cases}
\]
for which we can run the classical argument to conclude that the system exhibits exponential flocking with $\bv_i \to p(t^*)$. It also shows that the flock drifts to infinity to the right with the average velocity $p(t^*)$.  \\
\\
Suppose the momentum is positive for some time $t_0$. Then we can start time at $t_0$ and appeal to the above argument to prove that the system separates from the wall and exhibits flocking. We now show that if the momentum is always negative, the flock eventually settles to a state where:
\[
\bx_i (t) \to \bar{\bx}_i  \geq \ell.
\]
\\
Note that it does not suffice to show that $\bx_i(t) - \bx_j(t) \to \overline{\bx}_{ij}$. To prove that this happens, we start by showing that the velocity is integrable in time. 
\[\int_0^\infty |\bv|_\infty \dt < \infty.\] \\
\\
To prove that this happens, we start by showing that velocity (equivalently the momentum) is integrable in time. 
The derivative of the center of mass is, \[ \displaystyle \frac{1}{N} \sum_{i=1}^N \ddt \bx_i = \frac{1}{N} \sum_{i=1}^N \bv_i = p.\] \\
\\
Thus, and noting that $p = -|p|$, we have
\[\displaystyle \int_0^{\infty} |p(t)|dt = -\int_0^{\infty} p(t)dt = \lim_{t\to \infty} x.\] \\
\\
In particular, if velocity is not integrable, then $x\to -\infty$.\\
\\
If the average position of the agents goes to negative infinity, some agent must hit the wall, which is impossible.
\\
\\
Next, consider $|\vmin| = - \vmin$. As we have shown previously
\[
\ddt |\vmin|  \leq  - c_0 | \vmin| + c_0 | p|.
\]
Since $|p|$ is integrable we conclude that so is $|\vmin|$. 
Taking $\bv= \vmin$, we then have, \[\dot{\bv}\leq c_0 \bv+ c_0|p|.\] 
\\
Letting $f(s) = c_0(p(s))$, Duhamel gives us, 
\[ \bv(t) \leq e^{-c_0t}\bv_0+ \int_0^t e^{-c_0(t-s)}c_0(p(s))ds. \]
\\
Thus, 
\[ \int_0^\infty \bv(t) dt \leq \int e^{-c_0t}\bv_0dt + \int_0^\infty \int_0^t e^{c_0(t-s)}f(s)ds \thinspace dt.\] \\
By Fubini's theorem, we may rewrite the second integral as,
\[ \displaystyle \int_0^\infty f(s) \int_s^\infty e^{-c_0(t-s)} \thinspace dt \thinspace ds \]
where $s\leq t. $ \\
\\
Thus, $f$ is globally integrable. Knowing that $f$ is globally integrable, we can apply Duhamel in a similar fashion as above to show that $\vmax$ is globally integrable. This implies that all the $\bv_i$ are globally integrable. By the fundamental theorem of calculus, $\bx_i(t) \to \overline{\bx}_i(t)$ as t goes to infinity. If some $\overline{\bx}_i$ were less than $\ell$, the force would be bounded away from 0 for large enough time, which is impossible since force globally integrable. Thus, the $\bv_i$ go to 0 and the flock settles outside the interval $[0,\ell)$.
\\
\\
\section{Flocks confined inside a bounded interval}
We now consider the Cucker-Smale system on a bounded interval $[a,b]$. This 
case study is radically different from the previous one in  that the flock has 
no free directions to disperse. It results in more complex dynamics as the 
agents may reorganize again after bouncing off the walls. Our result in this case still states the ultimate alignment of the system, however storing flocking remains indeterminite and will be left to future research.

\begin{theorem}
	Any solution to the system \eqref{Cucker-Smale} with initial configuration $\bx_i(0) \in (a,b)$ will fulfill the following dynamic behavior:
	\begin{itemize}
		\item[(1)] It will never collide with the walls, 
		\[
		\inf_{i,t>0}\{ |\bx_i - a|, |\bx_i-b| \} > 0.
		\]
		\item[(2)] It will align
		\[
		A(t) = \max_{i,j} |\bv_i - \bv_j| \to 0.
		\]
	\end{itemize}
\end{theorem}

 To study this case 
we consider the full energy $\cE = \cK + \cP$ where $\cK$ and $\cP$ are defined 
as before, except the potential $U$ in this case is repulsive from both ends of 
the interval. The goal is to show that $\cE(t) \to 0$ as $t \to \infty$. As 
opposed the half-line case we cannot claim any rate in this case, in view of 
the complex interactions with the wall we remarked above. Since the 
communication kernel is non-degenerate $\phi >0$, on the compact set it implies 
that in fact there is a bound from below $\phi(r) \geq c_0$. From the energy equality
\[
\ddt \cE = - I_2 \leq - c_0 \cK
\]
we conclude that the kinetic energy is integrable in time :
\[
\int_0^\infty \cK(t) dt < \infty.
\]
Also, the energy itself $\cE$ is non-increasing. In particular, the potential energy remains bounded, $\cP \leq \cP_0$, from which we conclude as in the previous section that the agents do not crush into the walls
\[
\inf_{i,t>0}\{ |\bx_i - a|, |\bx_i-b| \} > \d.
\]
We also have globally bounded velocities
\[
\sup_{i, t>0} |\bv_i(t)| < M.
\]
To show that $\cK \to 0$ we will use the above integrability of the energy in time, together with an additional ingredient -- uniform continuity. In order to prove uniform continuity it suffices to show that $\cK$ stays Lipschits, in other words, $|\cK '| < \infty$. Indeed, if that's the case, and if $\cK \not \to 0$, then there is a sequence of times $\cK(t_i) > c_1$ for some $c_1 > 0$. Then there exists $\d>0$ so that $|\cK(t) | > c_1 /2$ for all $t\in [t_i - \d, t_i + \d]$. This implies that  
\[
\int_0^\infty \cK(t) dt \geq \sum_i \d c_1 = \infty, 
\]
a contradiction. Now, let us compute
\[
\cK' = - I_2 + \sum_i \bv_i \cdot \bF_i.
\]
Since $|I_2| \leq \cK |\phi|_\infty$, and by Cauchy-Schwarts,
\[
| \sum_i \bv_i \cdot \bF_i | \leq \sqrt{\cK} \bF < \infty,
\]
we arrive at the needed conclusion. The convergence $\cK \to 0$ already tells us that the flock aligns. Now we would like to show that all agents eventually settle outside the range of the potential $U$. We will do it, in fact, by showing that $U'(\bx_i) \to 0$ for each $i$. Since by our definition $U' = 0$ only outside the range of the potential, this will prove the desired result. In order to achieve the goal we introduce a new quantity -- work of force:
\[
W = \sum_i \bv_i U' (\bx_i).
\]
In view of the obtained information so far, we conclude that $|W|$ remains bounded for all time. Now, let us compute its derivative
\[
W' = \sum_{ij}\phi(\bx_i - \bx_j) [ \bv_i - \bv_j] \bv_i U' (\bx_i) - \sum_{i} |\bF_i|^2 + \sum_i U''(\bx_i) |\bv_i|^2.
\]
Note that $U''>0$, so the last term contributes as a non-negligible force.  However, since $|U''|< C$, in view of no collisions with the wall, we conclude that that term is integrable
\[
\int_0^\infty  \sum_i U''(\bx_i) |\bv_i|^2 dt \leq C \int_0^\infty \cK dt < \infty.
\]
For the first term, we apply the generalized Young inequality:
\[
\left| \sum_{ij}\phi(\bx_i - \bx_j) [ \bv_i - \bv_j] \bv_i U' (\bx_i) \right| \leq C_\e \cK + \e  \sum_{i} |\bF_i|^2.
\]
With $\e$ small we can see that the second part is absorbed by the forces the in the works-of-force budget law above, yet the first part is integrable.  In conclusion we obtain
\[
W ' \leq  - \sum_{i} |\bF_i|^2 + C \cK,
\]
which implies that 
\[
\int_0^\infty  \sum_{i} |\bF_i|^2 dt <\infty.
\]
Following the same set of ideas ws with the kinetic energy, we now only need to prove that $ |\bF_i|^2$ is uniformly continuous in order to conclude that $ |\bF_i| \to 0$. But,
\[
\ddt |\bF_i|^2 = U''(\bx_i)U'(\bx_i)\bv_i \leq C,
\]
in view of the non-collisions and boundedness of the velocities.

\section*{Acknowledgements} I would like to thank my faculty mentor, Professor Roman Shvydkoy, for his guidance and support. In addition, I would like to thank the University of Illinois at Chicago for offering the LAS Undergraduate Research Initiative, under which this project was undertaken.
\\
\\


\end{document}